\documentclass[a4paper,12pt]{article}
\usepackage{graphicx} 
\usepackage{psfrag}
\usepackage{amsfonts}
\usepackage{amssymb}
\usepackage{amsmath}
\usepackage{latexsym}
\usepackage[latin1]{inputenc}
\usepackage{ae}

\numberwithin{equation}{section}

\newcommand{\D}{{\displaystyle}}
\newcommand{\R}{{\mathbb R}}
\newcommand{\C}{{\mathbb C}}
\newcommand{\N}{{\mathbb N}}
\newcommand{\Z}{{\mathbb Z}}
\newcommand{\dokend}{\hfill \hbox{\vrule width 5pt height 5pt depth 0pt}}

\begin{document}
\thispagestyle{empty}
\author{ {\normalsize Matthias Kunik
\footnote{matthias.kunik@mathematik.uni-magdeburg.de}}\\
\small Institute for Analysis and Numerics, Otto-von-Guericke-Universit\"at\\
\small Postfach 4120 . D-39106 Magdeburg, Germany}
\title{Logarithmic Fourier integrals for the\\ Riemann Zeta Function}
\date{\today}
\maketitle



Preprint No. 01/2008,\\ 
Otto-von-Guericke University of Magdeburg,\\
Faculty of Mathematics.\\

\thispagestyle{empty}
\newpage
\setcounter{page}{1}
\begin{abstract}
\noindent
We use symmetric Poisson-Schwarz formulas
for analytic functions $f$ in the half-plane $\mbox{Re}(s)>\frac12$
with $\overline{f(\overline{s})}=f(s)$
in order to derive factorisation theorems for the Riemann zeta function.
We prove a variant of the Balazard-Saias-Yor theorem and obtain
explicit formulas for functions which are important
for the distribution of prime numbers. In contrast to 
Riemann's classical explicit formula, these
representations use integrals along the critical line
$\mbox{Re}(s)=\frac12$ and Blaschke zeta zeroes.

\noindent
\end{abstract}
\setcounter{section}{0}
\setcounter{equation}{0}
\section{Introduction}

Due to a theorem of Balazard-Saias-Yor \cite{BSY},
the Riemann-hypothesis is true if 
and only if the integral
\begin{align}\label{bsyintegral}
\Omega_{\zeta}:=\frac{1}{2\pi}\int \limits_{\mbox{Re}(w)=1/2} 
\frac{\log |\zeta(w)|}{|w|^2}\,|dw|
\end{align}
vanishes. These studies have their origin in Beurlings work
\cite{Beu}, see also \cite{BS}. Using that $(s-1)\zeta(s)/s^2$
belongs to the Hardy space $H_2(\mbox{Re}(s)>1/2)$
one can also show that the logarithmic integral \eqref{bsyintegral} 
is absolutely convergent, see Burnol \cite{Br1}, \cite{Br2}.
The general result about absolute convergence of the logarithmic integral
for Hardy spaces $H_p(\mbox{Im}(s)>0)$ is given for example
in the textbooks of Koosis \cite{Koo1} and Garnett \cite{Gar}.

The following factorisation formula presented in \cite{Br1} is valid 
for $\mbox{Re}(s)>1/2$
\begin{align}\label{burnfact0}
\frac{s-1}{s} \zeta(s)=
\zeta_B(s) \cdot B(s)\,,
\end{align}
with the function $\zeta_B$ in the right half plane $\mbox{Re}(s)>1/2$
given by
\begin{align}\label{burnfact1}
\zeta_B(s):=
\exp \left[
\frac{1}{2 \pi} \int \limits_{\mbox{Re}(w)=1/2}
\log \left| \zeta(w) \right| 
\frac{s+w-2sw}{s-w} \, \frac{|dw|}{|w|^2}
\right] \,,
\end{align}
and the Blaschke product
\begin{align}\label{burnfact2}
B(s):=
\prod \limits_{\substack{\rho \,:\, \zeta(\rho)=0\\ Re(\rho) > 1/2}}
\left\{
\frac{1-
\begin{displaystyle}\frac{s}{\rho}
\end{displaystyle}}{
\begin{displaystyle}
1-
\begin{displaystyle}
\frac{s}{1-\overline{\rho}}
\end{displaystyle}
\end{displaystyle}}^{~}
\cdot
\left|\frac{\rho}{1-\rho}  \right|
\right\}
\,.
\end{align}
It follows than by putting $s=1$ 
into \eqref{burnfact0}, \eqref{burnfact1}, \eqref{burnfact2} that
\begin{align}\label{burnfact3}
B(1)=
\prod \limits_{\substack{\rho \,:\, \zeta(\rho)=0\\ Re(\rho) > 1/2}}
\left| \frac{1-\rho}{\rho} \right|\,, \quad \Omega_{\zeta}=-\log B(1) \geq 0\,,
\end{align}
which gives the Balazard-Saias-Yor Theorem.

By using the symmetry of $\log |\zeta(1/2+iu)|$ we can also rewrite
$\zeta_B$ in the form
\begin{align}\label{zetabpoisson}
\zeta_B(s) =
\exp \left[
\frac{2}{\pi} \left(s-\frac12\right)
\int \limits_{0}^{\infty}
\frac{\log \left| \zeta(\frac12+iu) \right|} 
{u^2+(s-\frac12)^2} \, du\,
\right] \,, \quad \mbox{Re}(s)>1/2\,.
\end{align}
Thus the function $\zeta_B$ results from $\zeta$ by 
multiplication with $\frac{s-1}{s}$ and replacing
the Blaschke zeros $\rho$ of the zeta function with $\mbox{Re}(\rho)>\frac12$
by $1-\overline{\rho}$, 
i.e. by reflecting these zeros on the critical line
$\mbox{Re}(s)=1/2$ into the left half plane $\mbox{Re}(s)<1/2$.

Blaschke products and factorisation formulas 
play an important role in the theory of Hardy spaces, 
see for example the textbooks of Koosis \cite{Koo1}, 
Hoffmann \cite{Hoff} and Garnett \cite{Gar}. 
In the monographs \cite{Koo2,Koo3} of Koosis
the meaning of logarithmic integrals like \eqref{bsyintegral}
is highlighted by an important theoretical background
including results of Beurling, Malliavin and others, 
but the Riemann zeta function is not considered there.

In Section 2 we derive another symmetric Poisson-Schwarz formula
by using the Hadamard product decomposition of $(s-1) \zeta(s)$
and the Riemann-von-Mangoldt function $N(t)$ counting the nontrivial
zeta zeroes. The resulting formula is a counterpart of \eqref{zetabpoisson}. 
Then we prove a variant of the Balazard-Saias-Yor Theorem.

In Section 3 we apply the Fourier transform or Mellin's inversion formula on
\begin{equation}\label{logzetaa}
-\frac{\log \left((s-1)/s \, \zeta(s) \right)}{s(s-1)} \,, 
\quad \mbox{Re}(s)>1\,,
\end{equation}
and relate it to the Fourier transforms of
\begin{equation}\label{logzetab}
 \frac{\log |\zeta(\frac12+it)|}{t^2+\frac14}\,, \quad
 \frac{\pi N(t) - \vartheta(t)-2 \arctan(2t)}{t^2+\frac14}\,, \quad t \in \R\,,
\end{equation} 
by using the Blaschke zeta zeroes and a variant of Riemann's
explicite prime counting formula. If the Riemann-hypothesis is true,
then \eqref{logzetab} can be interpreted as the trace of the real- and
imaginary part of the analytical conti\-nu\-ation of
\eqref{logzetaa} on the critical line, respectively. 

Conversely, some important
number theoretic functions are obtained in terms of integrals
along the critical line by using the two expressions in \eqref{logzetab}.
   
\newpage

\section{A further factorisation of the Riemann zeta function
with a symmetric Poisson-Schwarz integral}
\setcounter{equation}{0}

In this section we first prove a Poisson-Schwarz representation which turns out to be
a dual counterpart of equations \eqref{burnfact0} and \eqref{zetabpoisson}.
In the next section we will combine these results 
with Fourier- or Mellin transform
techniques in order to obtain integrals along the critical line
for the representation of interesting explicite prime number formulas.

The desired Poisson-Schwarz formula results
from the Hadamard factorisation of $\zeta(s)$.
For this purpose we need the following theorem, 
which provides some information about the vertical distribution of the zeros 
of the zeta function in the critical strip. The following result 
was used by Riemann in \cite{Rm}, it was first shown 
by von Mangoldt \cite{Mn2} and then simplified by Backlund \cite{Bck1},
see the Edwards textbook \cite{Ed}\,.

{\bf Theorem (2.1)} {\em For $t \geq 0$ let $N(t)$ be the number of zeros
$\rho = \sigma +i\tau$ of the $\zeta$-function in the critical strip
$0 < \sigma <1$ with $0 \leq  \tau \leq t$, regarding the multiplicity
of the roots.  For $t < 0$ we put $N(t):=-N(-t)$, and
define the function $\vartheta: \R \to \R$ by
\begin{align}
\vartheta(t) := -\arctan(2t)-\frac{t}{2}(\gamma+\log \pi) +
\sum \limits_{k=1}^{\infty}
\left\{
\frac{t}{2k}- \arctan\frac{t}{2k+\frac12}
\right\}\,.\label{theta}
\end{align}

Then we have for $t \to \infty$ the asymptotic relations
\begin{equation}\label{nt1}
\begin{split}
N(t) = \frac{1}{\pi} \vartheta(t) + O(\log t) \,,\\
\sum \limits_{\substack{\rho \,:\, \zeta(\rho)=0\\ 0 < Re \, \rho < 1\\
|Im \, \rho| \leq t}}\,\frac{1}{|\rho|} = \frac{1}{2 \pi}\,(\log t)^2 
+ O(\log t)\,,\\
\sum \limits_{\substack{\rho \,:\, \zeta(\rho)=0\\ 0 < Re \, \rho < 1\\
|Im \, \rho| > t}}\,\frac{1}{|Im \, \rho|^2} 
=  \frac{1}{\pi}\,\frac{\log t}{t} + O(1/t)\,.
\end{split}
\end{equation}
}
\dokend

{\bf Remark:} The exact form of $\vartheta(t)$ is needed for later purposes
beside its asymptotic approximation 
\begin{align}\label{asymp}
\theta(t) \sim \frac{t}{2} \log \frac{t}{2 \pi}-
\frac{t}{2} -\frac{\pi}{8} + o(1)\,, \quad t \to \infty\,.
\end{align}

{\bf Theorem (2.2)} {\em We define for $\mbox{Re}(s)>1/2$
the analytic function $\zeta_C$ by
\begin{align}\label{zetacpoisson}
\zeta_C(s) :=
\exp \left[
s(s-1)\frac{2}{\pi} 
\int \limits_{0}^{\infty}
\frac{u \left( \pi N(u)-\vartheta(u)-2 \arctan(2u) \right) }
{(u^2+(s-\frac12)^2) \, (u^2+\frac14)}\, du\,
\right] \,,
\end{align}
as well as the symmetric counterpart $C(s)$ of the Blaschke product\,,
\begin{align}\label{cproduct}
C(s):=
\prod \limits_{\substack{\rho \,:\, \zeta(\rho)=0\\ Re(\rho) > 1/2}}
\frac{
\left(1-
\begin{displaystyle}
\frac{s}{\rho}
\end{displaystyle}
\right) \left(1-
\begin{displaystyle}
\frac{s}{1-\overline{\rho}}
\end{displaystyle}
\right)
}
{\begin{displaystyle}
\left(1 - 
\frac{2s}{1-{\overline{\rho}}+\rho}
\right)^{\large 2}
\end{displaystyle}
}
\,.
\end{align}
Then the following representation is valid for $\mbox{Re}(s)>1/2$\,,
\begin{align}\label{logzetaim}
\frac{s-1}{s}\zeta(s)=\zeta_C(s) \cdot C(s)\,.
\end{align}
}

{\bf Remarks:} Note that the integral in \eqref{zetacpoisson}
is well defined by the Riemann-von Mangoldt Theorem (2.1).
The function $\zeta_C$ results from $\zeta$ by 
multiplication with $\frac{s-1}{s}$ and projecting
the Blaschke zeros $\rho$ of the zeta function with $\mbox{Re}(\rho)>1/2$ 
on the critical line $\mbox{Re}(s)=1/2$. The Riemann-Hypothesis
is equivalent to $B(s)=C(s)\equiv 1$.

{\bf Proof:} We define the integration kernel
\begin{align}\label{KCkern}
K_C(s,u):=\frac{2}{\pi}\frac{s(s-1)\cdot u}
{(u^2+(s-\frac12)^2) \, (u^2+\frac14)}\,,
\end{align}
and for all $\alpha \in \R$ the characteristic function 
$\chi_{\alpha}:\R \to \R$ by

\begin{equation}  \label{randpuls}
\chi_{\alpha}(u) := \left\{ 
\begin{array}{ccl}
1 & , & u \geq \alpha \\ 
0 & , & u  < \alpha \,. \\ 
\end{array} \right.
\end{equation}

The following integrals can be obtained for $\mbox{Re}(s)>\frac12$
from elementary theory of the Poisson-Schwarz integral representations,

\begin{align}\label{quadfactor}
\exp \left[\pi
\int \limits_{0}^{\infty}
 \chi_{\alpha}(u) K_C(s,u)\, du\,
\right] 
=
\left(
1-\frac{s}{\frac12+i\alpha}
\right)
\left(
1-\frac{s}{\frac12-i\alpha}
\right)\,, \quad \alpha \geq 0\,,
\end{align}

\begin{align}\label{gammafactor}
\exp \left[
\int \limits_{0}^{\infty}
 \arctan\left( \frac{u}{a}\right) K_C(s,u)\, du\,
\right]
=1+\frac{s-1}{a+\frac12}\,, \quad \mbox{Re}(a) > 0\,,
\end{align}

\begin{align}\label{sminus1}
\int \limits_{0}^{\infty}
 u K_C(s,u)\, du = s-1\,.
\end{align}

Using \eqref{theta} and \eqref{gammafactor}, \eqref{sminus1} we obtain 

\begin{equation}\label{thetaKCur}
\begin{split}
\exp \left[
\int \limits_{0}^{\infty}
 \vartheta(u) K_C(s,u)\, du\,
\right] =
\frac{\pi^{-\frac{s-1}{2}}}{s}
\begin{displaystyle}
e^{-\gamma \frac{s-1}{2}}
\end{displaystyle}
\prod \limits_{k=1}^{\infty}
\frac{
\begin{displaystyle}
e^{\frac{s-1}{2k}}
\end{displaystyle}
}{1+\frac{s-1}{2k+1}}\,.
\end{split}
\end{equation}

Note that the integral on the left hand side
in \eqref{thetaKCur} is absolutely convergent
by the asymptotic relation \eqref{asymp}.

By using the general product representation
\begin{equation}\label{gammaprod}
\prod \limits_{k=1}^{\infty}
\left\{
\left( 1 + \frac{z}{k} \right) e^{-\frac{z}{k}}
\right\}=
\frac{e^{-\gamma z}}{\Gamma(z+1)}\,, \quad z \in \C\,,
\end{equation}
we can rewrite the right hand side of \eqref{thetaKCur} as
\begin{equation}\label{thet2}
\begin{split}
& \frac{
\pi^{-\frac{s-1}{2}}
\begin{displaystyle}
e^{-\gamma \frac{s-1}{2}}
\end{displaystyle}}{2^{s-1}}
\prod \limits_{k=0}^{\infty}
\frac{
\begin{displaystyle}
e^{\frac{s-1}{2k+1}}
\end{displaystyle}
}{1+\frac{s-1}{2k+1}} \nonumber\\
= &
\frac{
\pi^{-\frac{s-1}{2}}
\begin{displaystyle}
e^{-\gamma \frac{s-1}{2}}
\end{displaystyle}}{2^{s-1}}
\frac{e^{- \gamma \frac{s-1}{2}}}{\Gamma\left( \frac{s-1}{2} + 1 \right)}
\prod \limits_{k=1}^{\infty}
\frac{
\begin{displaystyle}
e^{\frac{s-1}{k}}
\end{displaystyle}
}{1+\frac{s-1}{k}} \nonumber \\
= & \frac{\pi^{-\frac{s-1}{2}} e^{-\gamma (s-1)}}{2^{s-1} 
\Gamma \left( \frac{s+1}{2}\right)}
\frac{\Gamma(s)}{e^{-\gamma (s-1)}}
\nonumber \\
= & \Gamma\left(\frac{s}{2} \right)\, \pi^{-\frac{s}{2}}\,,
\end{split}
\end{equation}
where the last step follows from the duplication formula 
for the $\Gamma$ function, see for example the textbook 
of Andrews et al. \cite{AAR}.

We obtain from \eqref{thetaKCur} for $\mbox{Re}(s)> \frac12$ that
\begin{equation}\label{thetaKC}
\exp \left[
\int \limits_{0}^{\infty}
 \vartheta(u) K_C(s,u)\, du\,
\right]
= \Gamma\left(\frac{s}{2} \right)\, \pi^{-\frac{s}{2}}\,.
\end{equation}
From \eqref{gammafactor} it also follows with $a=\frac12$ that
\begin{equation}\label{atanKC}
\exp \left[
\int \limits_{0}^{\infty}
2 \arctan(2u) K_C(s,u)\, du\,
\right]
= s^2\,, \quad \mbox{Re}(s) > \frac12\,.
\end{equation}
If we denote the positive imaginary parts of the nontrivial zeta zeros by
\begin{equation}\label{roots1}
0 < t_1 \leq t_2 \leq t_3 \leq \cdots\,,
\end{equation}
where the imaginary parts are listed according
to the multiplicity of each root, 
then we can rewrite $N(u)$ for $u>0$ as
\begin{equation}\label{roots2}
N(u) = \sum \limits_{n=1}^{\infty} \chi_{t_n}(u)\,,
\end{equation}
and obtain from \eqref{quadfactor} that
\begin{equation}\label{NKC}
\exp \left[
\int \limits_{0}^{\infty}
 \pi N(u) K_C(s,u)\, du\,
\right]
= \prod \limits_{n=1}^{\infty}
\left\{
\left(
1-\frac{s}{\frac12+it_n}
\right)
\left(
1-\frac{s}{\frac12-it_n}
\right)
\right\}\,.
\end{equation}
Now we use the definition of the analytic functions $\zeta_C$
in \eqref{zetacpoisson} and $C$ in \eqref{cproduct}
and calculate from \eqref{NKC}, \eqref{thetaKC} and \eqref{atanKC}
for $\mbox{Re}(s)>\frac12$ the product
\begin{equation}\label{zetaccprod}
\begin{split}
\zeta_C(s)\,C(s) & = 
\frac{
\prod \limits_{\substack{\rho \,:\, \zeta(\rho)=0\\ 
Re(\rho), Im(\rho) > 0}}
\left\{
\left(
1-\frac{s}{\rho}
\right)
\left(
1-\frac{s}{1-\rho}
\right)
\right\}}
{\Gamma\left(\frac{s}{2} \right)\, \pi^{-\frac{s}{2}}\, s^2}\,.
\end{split}
\end{equation}
The Theorem results
from \eqref{zetaccprod} with Hadamard's product decomposition 
of Riemann's function $\xi : \C \to \C$ with
\begin{equation}\label{xihadamard}
\xi(s) := \frac{s(s-1)}{2}\,\pi^{-\frac{s}{2}}\,\Gamma(\frac{s}{2})\,\zeta(s)
= \frac12 \, \lim \limits_{T \to \infty} 
\prod \limits_{\substack{\rho \,:\, \zeta(\rho)=0\\ 0 < Re \, \rho < 1\\
|Im \, \rho| \leq T}}\,\left( 1-\frac{s}{\rho} \right)\,.
\end{equation}
\dokend

{\bf Corollary (2.3)} ~ {\em With the integration kernel \eqref{KCkern} we have
\begin{align}\label{xipoisson}
\xi(s) = \frac{C(s)}{2} \cdot \exp \left[
\int \limits_{0}^{\infty}
\pi N(u) K_C(s,u)\, du\,
\right] \,, \quad \mbox{Re}(s)>\frac12\,.
\end{align}}

{\bf Remark:}~ For the proof of Theorem (2.2) we have directly used
the Riemann von Mangoldt Theorem (2.1) and the Hadamard products
for $\zeta$ and $\xi$. In contrast,
the integral \eqref{bsyintegral}
is divergent if we replace there $\zeta(w)$ 
by $\Gamma(w/2)$ or by Riemann's function $\xi(w)$.
In this case even the Cauchy limit
$$
\lim \limits_{T \to \infty} \int \limits_{-T}^{T}
\frac{\log |\xi(1/2+it)|}{1/4+t^2}\,dt
$$
is divergent due to the asymptotic behaviour
$$
\lim \limits_{|t| \to \infty} 
\frac{\log |\Gamma(1/4+\frac{it}{2})|}{\frac{\pi}{4}|t|}=-1
$$
of the Gamma function. Thus it is not possible
to evaluate \eqref{bsyintegral} or \eqref{zetabpoisson}
by applying Hadamard's product representation 
for $\zeta(w)$ directly on $\log |\zeta(w)|$.

We have obtained Theorem (2.2) completely independent from
the theory of Hardy spaces and the results of Balazard-Saias-Yor.
Note that for all $t \in \R$
\begin{align}\label{spur1}
\zeta(\frac12 +it)=\exp \left[ \log |\zeta(\frac12 +it)|+
i \left(\pi N(t) - \vartheta(t)- \pi \, \mbox{sign}(t)\right) \right]\,.
\end{align}
Moreover, if the Riemann-hypothesis is true, then an analytical logarithm
of $\frac{s-1}{s} \zeta(s)$ is defined for $\mbox{Re}(s)>\frac12$
with real values for $s>1/2$ and trace 
\begin{align}\label{spur2}
\log |\zeta(\frac12 +it)|+
i \left(\pi N(t) - \vartheta(t)-2 \arctan(2t)\right)
\end{align}
on $s=\frac12+it$.
Now we obtain from Theorem (2.2) the following 
counterpart of the Balazard-Saias-Yor theorem, namely

{\bf Theorem (2.4)} ~{\em We have the two inequalities
\begin{align}\label{BSY1counter}
-\frac{1}{\pi} 
\int \limits_{0}^{\infty}
\frac{\log \left| \zeta(\frac12 + iu) \right|}
{(u^2+\frac14)^2}\, du \leq \gamma-1\,,
\end{align}
\begin{align}\label{BSY2counter}
\frac{2}{\pi} 
\int \limits_{0}^{\infty}
\frac{u \left( \pi N(u)-\vartheta(u)-2 \arctan(2u) \right) }
{(u^2+\frac14)^2}\, du \geq \gamma -1\,.
\end{align}
In both cases equality holds if and only if
the Riemann hypothesis is valid.}

{\bf Proof:} If we denote the left hand side
of the inequalities \eqref{BSY1counter}, \eqref{BSY2counter} 
by $J_1$ and $J_2$, respectively,
and put $s=1$ in the logarithmic derivatives
of \eqref{burnfact0} and \eqref{logzetaim}, 
then we obtain
\begin{align}\gamma-1=
\lim \limits_{s \to 1}
\left[
\frac{1}{s-1}-\frac{1}{s}+\frac{\zeta'(s)}{\zeta(s)}
\right] = J_1+2\Omega_{\zeta} +\frac{B'(1)}{B(1)}
= J_2+\frac{C'(1)}{C(1)}\,.
\end{align}
We will first prove that $2\Omega_{\zeta}+\frac{B'(1)}{B(1)} \geq 0$,
and that $2\Omega_{\zeta}+\frac{B'(1)}{B(1)}=0$ is equivalent to the Riemann
hypothesis. Put
\begin{align}\label{crho}
f_{\rho}:=
2 \log \left| \frac{\rho}{1-\rho} \right| + \frac{1}{1-\rho}- 
\frac{1}{\overline{\rho}}
\end{align}
for any $\rho = \sigma + i \tau$ with 
$\frac12 < \sigma \leq 1$ and $\tau \in \R$ with $|\tau| >1$.
Then we obtain
\begin{align}\label{frhoeval}
f_{\rho}+f_{\overline{\rho}}= 2\int \limits_{\frac12}^{\sigma}
\frac{x^2(1-x)^2+\tau^2(\tau^2-6x^2+6x-1)}
{(x^2+\tau^2)^2 \left((1-x)^2+\tau^2\right)^2}
\,dx > 0\,,
\end{align}
such that the first part of the Theorem results from
\begin{align}\label{frhoeval2}
\frac{B'(1)}{B(1)}+2\Omega_{\zeta}=
\sum \limits_{\substack{\rho \,:\, \zeta(\rho)=0\\ Re(\rho) > 1/2}} f_{\rho}\,.
\end{align}
It remains to show that $C'(1)/C(1)\leq 0$ and that
the equality $C'(1)/C(1)= 0$ holds if and only if 
the Riemann hypothesis is valid. Put
\begin{align}\label{crho}
C_{\rho}(s):=
\frac{
\left(1-
\begin{displaystyle}
\frac{s}{\rho}
\end{displaystyle}
\right) \left(1-
\begin{displaystyle}
\frac{s}{1-\overline{\rho}}
\end{displaystyle}
\right)
}
{\begin{displaystyle}
\left(1 - 
\frac{2s}{1-{\overline{\rho}}+\rho}
\right)^{\large 2}
\end{displaystyle}
}
\end{align}
for any $\rho = \sigma + i \tau$ with 
$\frac12 < \sigma \leq 1$ and $\tau \in \R$.
Then we obtain
\begin{equation}\label{logcrho}
\frac{C'_{\rho}(1)}{C_{\rho}(1)}+
\frac{C'_{\overline{\rho}}(1)}{C_{\overline{\rho}}(1)}=
-\frac{2\left(3\tau^2-\sigma(1-\sigma)\right) \left(\sigma-\frac12\right)^2}
{((1-\sigma)^2+\tau^2)(\sigma^2+\tau^2)(\frac14+\tau^2)}\,.
\end{equation}
For $\begin{displaystyle}|\tau|>\frac{1}{\sqrt{12}}\end{displaystyle}$ 
and $\frac12<\sigma \leq 1$ this expression is always negative.
Thus we have shown the Theorem. \dokend\\

{\bf Theorem (2.5)} ~ {\em For $t \geq 0$ 
except on a discrete singular set we define the number $N_B(t)$
of Blaschke zeta zeros $\rho$ with $\mbox{Re}(\rho)>\frac12$
and $|\mbox{Im}(\rho)| \leq t$ as well as the quantities
\begin{align}\label{N1}
N_1(t):=\frac{1}{\pi}
\left[
\vartheta(t)+2 \arctan(2t)
\right]\,,
\end{align}
\begin{align}\label{N2}
N_2(t):=-\frac{1}{2\pi^2}
\frac{d}{dt}\,\int \limits_{-\infty}^{\infty}
\log \left| 1- \frac{t^2}{u^2} \right| \,
\log \left| \zeta(\frac12+iu) \right|\, du\,,
\end{align}
\begin{align}\label{N3}
N_3(t):=\frac{1}{2\pi}
\int \limits_{-t}^{t}
\frac{B'(\frac12+iu)}{B(\frac12+iu)}\, du\,.
\end{align}
If we extend $N_B$ and $N_1, N_2, N_3$ as odd functions to the whole real
axis, then we obtain for almost all $t \in \R$ that
\begin{align}\label{Nsumme}
N(t) - N_B(t) = N_1(t) + N_2(t) + N_3(t)\,. 
\end{align}}
{\bf Remarks:} 

(a) Without using any information about the horizontal distribution of the
zeta zeroes in the critical strip $0 < \mbox{Re}(s)<1$, one can
employ Theorem (2.1) and estimations for the logarithmic derivative of
the Blaschke product on the critical line
to prove the following asymptotic relations for $t \to \infty$,
\begin{align}\label{NBN3}
N_2(t) = O(\log^2 t)\,, \quad
N_3(t) + N_B(t) = O(\log^2 t)\,. 
\end{align}
(b) If the Riemann hypothesis is valid, then we obtain
from Theorem (2.1) with $N_3=N_B \equiv 0$
the better result $N_2(t) = O(\log t)$ for $t \to \infty$.

{\bf Proof:} For $s \in \C \setminus \{1\}$ with $\zeta(s) \neq 0$ we 
obtain from \eqref{xihadamard} 
\begin{align}\label{N49a}
\frac{\xi'(s)}{\xi(s)}=\frac{1}{s-1}-\frac12(\gamma+\log \pi)
+\sum \limits_{n=1}^{\infty} \left(\frac{1}{2n}-\frac{1}{s+2n}\right)
+ \frac{\zeta'(s)}{\zeta(s)}\,, 
\end{align}
and note the functional equation 
\begin{align}\label{N49d}
\frac{\xi'(1-s)}{\xi(1-s)}= -\frac{\xi'(s)}{\xi(s)}\,. 
\end{align}
The logarithmic derivative of \eqref{burnfact0}, \eqref{zetabpoisson} is
given for $\mbox{Re}(s)>\frac12$ with $s \neq 1$ and $\zeta(s) \neq 0$ by
\begin{align}\label{N49b}
\frac{\zeta'(s)}{\zeta(s)}=\frac{1}{s}-\frac{1}{s-1}
-\frac{1}{\pi}
\int\limits_{-\infty}^{\infty}
\frac{\log \left|\zeta(\frac12\pm iu)\right|}{(s-\frac12-iu)^2}\,du
+\frac{B'(s)}{B(s)}\,. 
\end{align}
Thus we obtain for $\mbox{Re}(s)>\frac12$ from \eqref{N49a} that
\begin{align}\label{N49c}
\frac{\xi'(s)}{\xi(s)}&=\frac{1}{s}-\frac12(\gamma+\log \pi)
+\sum \limits_{n=1}^{\infty} \left(\frac{1}{2n}-\frac{1}{s+2n}\right)
\nonumber\\
&-\frac{1}{\pi}
\int\limits_{-\infty}^{\infty}
\frac{\log \left|\zeta(\frac12\pm iu)\right|}{(s-\frac12-iu)^2}\,du
+ \frac{B'(s)}{B(s)}\,. 
\end{align}
For $\eta,t >0$ we define the rectangular positive oriented integration path 
$\it{P}_{\eta,t}$ consisting on the straight line segments
$[\frac12+\eta-it,\frac12+\eta+it]$,
$[\frac12+\eta+it,\frac12-\eta+it]$,
$[\frac12-\eta+it,\frac12-\eta-it]$,
$[\frac12-\eta-it,\frac12+\eta-it]$.
It is centred around $s=\frac12$, is the boundary of the rectangle
$|\mbox{Re}(s)-\frac12| \leq \eta$, $|\mbox{Im}(s)| \leq t$
and contains the critical strip for any $\eta > \frac12$.
We decompose $\it{P}_{\eta,t}$ into two parts $\it{P}_{\eta,t}^{\pm}$
for the complex numbers $s$ with $\mbox{Re}(s)\geq\frac12$ and
$\mbox{Re}(s)\leq\frac12$ respectively,
\begin{align}\label{pathpm}
{\it P}_{\eta,t}^{+}:=[\frac12-it,\frac12+\eta-it] \oplus
[\frac12+\eta-it,\frac12+\eta+it]\oplus
[\frac12+\eta+it,\frac12+it]\,,\nonumber\\
{\it P}_{\eta,t}^{-}:=[\frac12+it,\frac12-\eta+it] \oplus
[\frac12-\eta+it,\frac12-\eta-it] \oplus
[\frac12-\eta-it,\frac12-it]\,.
\end{align}
Let be $t>0$ given with $\zeta(\frac12+it)\neq0$. Since $\zeta$ and $\xi$
have only isolated zeros, we can choose $0 < \eta < \frac12$ small enough
such that the compact rectangle
$|\mbox{Re}(s)-\frac12| \leq \eta$, $|\mbox{Im}(s)| \leq t$
contains only zeros of the $\zeta$-function on the critical line.
Then we obtain due to \eqref{xihadamard} and \eqref{N49d} that
\begin{align}\label{N49g}
\frac{1}{2\pi i} \int \limits_{\it{P}_{\eta,t}^{\pm}}
\frac{\xi'(s)}{\xi(s)}\,ds
=N(t)-N_B(t)\,.
\end{align}
The direct integration of \eqref{N49c}
is proble\-matic due to the presence of the pole singularities
from the $\zeta$-zeros on the critical line. Instead of this 
we evaluate the primitive $\int \limits_0^t ( N(\tau)-N_B(\tau) )\,d\tau$.
We have assumed $\zeta(\frac12 \pm it) \neq 0$,
and the integral on the left hand side in \eqref{N49g}
does not depend on $\eta$ for sufficiently small $\eta >0$, and therefore
\begin{align}\label{N49l}
\int \limits_0^t ( N(\tau)-N_B(\tau) )\,d\tau
=\frac{1}{2 \pi}\,\lim \limits_{\eta \to 0^+}
\int \limits_{0}^{t} \int \limits_{-\tau}^{\tau}
\frac{\xi'(\frac12+\eta+i\vartheta)}{\xi(\frac12+\eta+i\vartheta)}
\,d \vartheta \,d\tau
\end{align}
by using for the short integration paths
$[\frac12-it,\frac12+\eta-it]$ and $[\frac12+\eta+it,\frac12+it]$ of 
${\it P}_{\eta,t}^{+}$ the well known standard estimates.

First we regard for each $u \in \R$ the double integral
\begin{equation}\label{doppelint}
\begin{split}
&\frac{1}{2 \pi} \int \limits_{0}^{t}
\int \limits_{-\tau}^{\tau}
\frac{d \vartheta}{(\eta+i(\vartheta-u))^2} d \tau \\
&=\frac{1}{2 \pi}
\left[\,
\log(\eta+i(t-u))+\log(\eta-i(t+u))-2\log(\eta-iu)\,
\right]\,.
\end{split}
\end{equation}
Its real part is an even function on $u \in \R$, and its imaginary part
an odd function on $u$. In the limit $\eta \to 0^+$ we obtain for $u \neq 0$
and $|u| \neq t$ that
\begin{equation}\label{doppelintlim}
\begin{split}
&\lim \limits_{\eta \to 0^+} \frac{1}{2 \pi}
\left[\,
\log(\eta+i(t-u))+\log(\eta-i(t+u))-2\log(\eta-iu)\,
\right]\\
&=\frac{1}{2 \pi} \log \left|1-\frac{t^2}{u^2} \right| + \frac{i}{4}
\left[
\mbox{sign}(t-u)-\mbox{sign}(t+u)+ 2 \, \mbox{sign}(u)
\right]
\,.
\end{split}
\end{equation}
Moreover, we have
\begin{align}\label{N49k}
\frac{1}{2\pi i} \int \limits_{0}^{t}
\int \limits_{\frac12+\eta-i\tau}^{\frac12+\eta+i\tau}
\,ds\,d\tau
=\frac {t^2}{2 \pi}\,,
\end{align}
and for all $\alpha \geq 0$
\begin{align}\label{N49j}
\frac{1}{2\pi i} \int \limits_{0}^{t}
\int \limits_{\frac12+\eta-i\tau}^{\frac12+\eta+i\tau}
\frac{ds}{s+\alpha}\,d\tau 
=\frac {1}{\pi}\,
\left[
t \arctan\frac{t}{\frac12+\alpha+\eta}\right.  \nonumber \\
\left. -\frac{\frac12+\alpha+\eta}{2}
\log \left(1+\frac{t^2}{(\frac12+\alpha+\eta)^2} \right)
\right]\,.
\end{align}
Thus we obtain by direct calculation from the last two equations 
and \eqref{theta} 
\begin{equation}\label{N1form}
\begin{split}
&\lim \limits_{\eta \to 0^+}\frac{1}{2\pi i} \int \limits_{0}^{t}
\int \limits_{\frac12+\eta-i\tau}^{\frac12+\eta+i\tau}
\left[
\frac{1}{s}-\frac12(\gamma+\log \pi)
+\sum \limits_{n=1}^{\infty} \left(\frac{1}{2n}-\frac{1}{s+2n}\right)
\right]
\,ds\,d\tau\\
&= \int \limits_0^t N_1(\tau)\,d \tau\,.
\end{split}
\end{equation}
Next we use \eqref{doppelint}, \eqref{doppelintlim} and obtain 
from Fubini's theorem and the Lebesgue dominated convergence theorem that
\begin{equation}\label{N2form}
\begin{split}
\lim \limits_{\eta \to 0^+}\frac{1}{2\pi i} \int \limits_{0}^{t}
\int \limits_{\frac12+\eta-i\tau}^{\frac12+\eta+i\tau}
\left[-\frac{1}{\pi}
\int\limits_{-\infty}^{\infty}
\frac{\log \left|\zeta(\frac12\pm iu)\right|}{(s-\frac12-iu)^2}\,du
\right]
\,ds\,d\tau
= \int \limits_0^t N_2(\tau)\,d \tau\,.
\end{split}
\end{equation}
Finally we recall that for $\eta > 0$ sufficiently small
the rectangle bounded by $\it{P}_{\eta,t}^{\pm}$ does 
not contain Blaschke $\zeta$-zeros, such that
\begin{equation}\label{N3form}
\begin{split}
\lim \limits_{\eta \to 0^+}\frac{1}{2\pi i} \int \limits_{0}^{t}
\int \limits_{\frac12+\eta-i\tau}^{\frac12+\eta+i\tau}
\frac{B'(s)}{B(s)}
\,ds\,d\tau
= \int \limits_0^t N_3(\tau)\,d \tau\,.
\end{split}
\end{equation}
From \eqref{N49c}, \eqref{N49l} and 
\eqref{N1form}, \eqref{N2form}, \eqref{N3form}
we conclude \eqref{Nsumme}. \dokend

\section{Logarithmic Fourier integrals and their relation
to the distribution of prime numbers}
\setcounter{equation}{0}

In this section we derive Fourier-Mellin transforms
from the boundary integral formulas obtained in Section 2
which give interesting relations to the distribution of prime numbers.

As a first basic building block we need a representation theorem 
for the Mellin transforms of certain functions involving the
exponential integral, which is generally useful 
for the study of Hadamard's product decomposition 
of entire functions with appropriate growth conditions. 

One of the various representations of the exponential integral is
\begin{align}\label{intro1}
{\rm Ei}(z) := \gamma + \log z + {\rm Ei}_0(z)
\end{align}
with Euler's constant 
$\begin{displaystyle}
\gamma = \lim \limits_{n \to \infty} 
\left(
\sum \limits_{k=1}^n \frac{1}{k} - \log n
\right)
\end{displaystyle}$ 
and the entire function
\begin{align}\label{intro2}
{\rm Ei}_0(z) := \sum \limits_{k=1}^{\infty}\frac{z^k}{k \cdot k!} = 
\int \limits_{0}^{z} \frac{e^t - 1}{t}\,dt =
\int \limits_{0}^{1} \frac{e^{u z} - 1}{u}\,du\,.
\end{align}
In contrast to the entire function ${\rm Ei}_0$, the functions
${\rm Ei}$ and $\log$ are only defined on the cut plane
$\begin{displaystyle}
\C_{-} := \C \setminus (-\infty,0]\,.
\end{displaystyle}$ 

The representation theorem formulated below 
was proved in \cite{Ku1}, which also contains
its application to Riemann's explicit formula for the number of primes less 
than a given limit $x > 1$ and related formulas.\\

{\bf Theorem (3.1)} {\em Let $s$ be any complex number with ${\rm Re}(s)>0$.
\begin{itemize}
\item[(a)] The expression 
$ 
(\gamma + \log(\log x))/x^{s+1}  
$
is Lebesgue integrable on the interval $(1,\infty)$, 
and there holds the relation
$$
\frac{1}{s} = \exp \left( s\,
\int \limits_1^{\infty} \frac{\gamma + \log(\log x)}{x^{s+1}}\,dx \right)\,.
$$
\item[(b)]
The expression
$
{\rm Ei}(\log x)/x^{s+1}  
$
is Lebesgue integrable on the interval $(1,\infty)$, 
and with ${\rm Li}(x):={\rm Ei}(\log x)$ there holds the relation
\begin{align*}
\frac{1}{s-1} = \exp \left( s\,
\int \limits_1^{\infty} \frac{{\rm Li}(x)}{x^{s+1}}\,dx \right)\,.
\end{align*}
\item[(c)] Assume that $\rho \in \C \setminus [0,\infty)$
and ${\rm Re}(\rho) < {\rm Re}(s)$\,. Then 
\begin{align*}
1-\frac{s}{\rho} = \exp \left( - s\,
\int \limits_1^{\infty} 
\frac{\gamma + \log(\log x) + \log(-\rho) + {\rm Ei}_0({\rho} \log x)}
{x^{s+1}}\,dx \right)\,.
\end{align*} \dokend
\end{itemize}}
If we note that for $u>0$ and $\mbox{Re}(s)>\frac12$
\begin{equation}\label{sinlogint}
\begin{split}
\int \limits_{1}^{\infty} \frac{\sqrt{x} \, \sin(u \log x)}{x^{s+1}}\,dx
=\frac{u}{u^2+(s-\frac12)^2}
\end{split}
\end{equation} 
and define $f_C : (1,\infty) \to \R$ by
\begin{equation}\label{fC}
\begin{split}
f_C(x):=
\frac{2 \sqrt{x}}{\pi}\int \limits_{0}^{\infty} 
\frac{ \pi N(u) - \vartheta(u)-2 \arctan(2u)  
}{u^2+\frac14} \cdot \sin(u \log x) \, du \,,
\end{split}
\end{equation} 
then we obtain from \eqref{nt1} and \eqref{zetacpoisson} with Fubini's theorem
\begin{equation}\label{sincosint}
\begin{split}
\zeta_C(s)=
\exp\left[ s(s-1)\int \limits_{1}^{\infty} \frac{f_C(x)}{x^{s+1}}\,dx  \right]\,.
\end{split}
\end{equation} 
Next we define for $-\alpha \notin \C_{-}$
the functions $\varphi_{\alpha}\,, \Phi_{\alpha}  : (1, \infty) \to \C$ by
\begin{equation}\label{phialpha}
\begin{split}
\varphi_{\alpha}(x) &:= \gamma + \log(\log x) + \log(-\alpha)
+\mbox{Ei}_0(\alpha \log x)\,,\\
\Phi_{\alpha}(x) &:= x \int \limits_{1}^{x} \frac{\varphi_{\alpha}(y)}{y^2}dy\,.
\end{split}
\end{equation} 
From Theorem (3.1) above we obtain by partial integration 
that for $\mbox{Re}(s)>1$ and $\mbox{Re}(s) > \mbox{Re}(\alpha) > 0$ 
\begin{equation}\label{Phimellin}
1-\frac{s}{\alpha}= \exp \left[
 -s(s-1) \int \limits_{1}^{\infty} \frac{\Phi_{\alpha}(x)}{x^{s+1}} \, dx\right]\,.
\end{equation} 
The integral expression for $\Phi_{\alpha}$ in \eqref{phialpha}
can be solved explicitely. For any $\alpha \in \C$ with 
$\mbox{Im}(\alpha) \neq 0$ and for all $x>1$ we obtain that
\begin{equation}\label{phialphamod}
\begin{split}
\Phi_{\alpha}(x) &= x \cdot \varphi_{\alpha -1}(x)-\varphi_{\alpha}(x) + 
x \cdot \left[\,\log\left(-\alpha\right)-\log\left(-(\alpha-1)\right)\,\right] \,.
\end{split}
\end{equation} 
Equation \eqref{phialphamod} can be checked easily by forming the
derivative of $\Phi_{\alpha}(x)/x$ with respect to $x$
and by regarding that $\D \lim \limits_{x \to 1}\Phi_{\alpha}(x)=0$.
We have not expressed ${\Phi}_{\alpha}(x)$
in terms of the exponential integral, because $\varphi_{\alpha}(x)$
in \eqref{phialpha} has a better asymptotic behaviour than 
$\mbox{Ei}(\alpha \log x)$ for $\mbox{Im}(\alpha) \neq 0$ and $x>1$,
\begin{equation}\label{phiasym1}
\begin{split}
\varphi_{\alpha}(x) &= \mbox{Ei}(\alpha \log x) - 
i \, \pi \, \mbox{sign}(\mbox{Im}(\alpha))\\
&= \frac{x^{\alpha}}{\alpha \log x}+
x^{\alpha} \int \limits_0^{\infty}
\frac{e^{- y}}{(\alpha \log x - y)^2}\,dy\,,\\
& \left|\,  x^{\alpha} \int \limits_0^{\infty}
\frac{e^{- y}}{(\alpha \log x - y)^2}\,dy \,\right| 
\leq \frac{x^{\mbox{Re}(\alpha)}}{|\mbox{Im}(\alpha)|^2 \, {\log}^2 x}\,.
\end{split}
\end{equation}
A proof of the following theorem can be found in the textbook \cite{Ed}
of Edwards.\\

{\bf Theorem (3.2)}~{\em For any integer number $n \geq 1$  we define the 
von Mangoldt function
\[
\Lambda(n):=\left\{
\begin{array}{c@{\quad,\quad}l}
\ln p & \mbox{for}~ n=p^m, ~m \geq 1, ~ p ~ \mbox{prime}\\
0 & \mbox{otherwise}\,,
\end{array}
\right.
\]
and for $x \geq 1$ the functions
\[
\psi(x) := \sum \limits_{n \leq x} \Lambda(n)\,, \quad
\pi_*(x) := 
\sum \limits _{1 < n \leq x}\frac{\Lambda(n)}{\ln n} =
\sum \limits_{n=1}^{\infty}\frac{\pi(\sqrt[n]{x})}{n} \,,
\]
where $\pi(x)$ is the number of primes $\leq x$.
Then we obtain for ${\rm Re}(s)>1$
\begin{align}
\zeta(s) &= \exp(\,\sum \limits _{n=2}^{\infty}\frac{\Lambda(n)}
{\ln n}n^{-s}\,) 
= \exp(\, s\int \limits_{1}^{\infty}\frac{\pi_*(x)}{x^{s+1}}\,dx \,)
\label{zetapi}\\
-\frac{\zeta'(s)}{\zeta(s)} &= 
\sum \limits _{n=2}^{\infty}\Lambda(n) \, n^{-s} = 
s\,\int \limits_{1}^{\infty}\frac{\psi(x)}{x^{s+1}}\,dx\,. \label{zetapsi}
\end{align}} \dokend

From Theorem (3.1)(a,b) and \eqref{zetapi} we obtain with
the entire function $\mbox{Ei}_0$ defined in \eqref{intro2} 
for $\mbox{Re}(s)>1$ that
\begin{equation}\label{zetap1}
\frac{s-1}{s}\,\zeta(s)=
\exp\left[\, s\int \limits_{1}^{\infty}
\frac{\pi_*(x)-\mbox{Ei}_0(\log x)}{x^{s+1}}\,dx \,
\right]\,.
\end{equation}
If we apply partial integration on the last integral,
we can also rewrite \eqref{zetap1} in the form
\begin{equation}\label{zetap2}
\frac{s-1}{s}\,\zeta(s)=
\exp\left[\, s(s-1) \int \limits_{1}^{\infty}\frac{f_{*}(x)}{x^{s+1}}\,dx \,
\right]\,,
\end{equation}
where the function $f_{*}: (1, \infty) \to \R$ is given by
\begin{equation}\label{zetap3}
\begin{split}
f_{*}(x)&=x\,\int \limits_{1}^{x} \frac{\pi_*(y)-\mbox{Ei}_0(\log y)}{y^2}\,dy\\
&=x\,\left(
\sum \limits_{n \leq x}
\frac{\Lambda(n)}{n \log n}+\mbox{Ei}_0(-\log x)\right)-
\left( \pi_*(x)-\mbox{Ei}_0(\log x) \right)
\,.
\end{split}
\end{equation}

Now we are able to prove the following result,

{\bf Theorem (3.3)}~{\em For $\alpha \in \C$ 
with $\mbox{Im}(\alpha) \neq 0$ we define 
$\tilde{\Phi}_{\alpha}  : (1, \infty) \to \C$ by
\begin{equation}\label{tildephialpha}
\tilde{\Phi}_{\alpha}(x) := x \, \varphi_{\alpha-1}(x)-\varphi_{\alpha}(x)\,.
\end{equation} 
Then we obtain for all $x>1$ for the function $f_*$ in \eqref{zetap3}
\begin{itemize}
\item[(a)]
\begin{equation}\label{fsternsin}
\begin{split}
f_*(x) & =\frac{2 \sqrt{x}}{\pi}\int \limits_{0}^{\infty} 
\frac{ \pi N(u) - \vartheta(u)-2 \arctan(2u)}{u^2+\frac14} 
\cdot \sin(u \log x) \, du\\
&-\sum \limits_{\substack{\rho \,:\, \zeta(\rho)=0\\ \mbox{Re}(\rho) > \frac12}}
\left[
\tilde{\Phi}_{\rho}(x)+\tilde{\Phi}_{1-\overline{\rho}}(x)
-2\tilde{\Phi}_{\frac12(1+\rho-\overline{\rho}\,)}(x)
\right]\,,
\end{split}
\end{equation}
\item[(b)]
\begin{equation}\label{fsterncos}
\begin{split}
f_*(x)  &=-\frac{2 \sqrt{x}}{\pi}\int \limits_{0}^{\infty} 
\frac{ \log |\zeta(\frac12+iu)|}{u^2+\frac14} 
\cdot \cos(u \log x) \, du\\
&-\sum \limits_{\substack{\rho \,:\, \zeta(\rho)=0\\ \mbox{Re}(\rho) > \frac12}}
\left[
\tilde{\Phi}_{\rho}(x)-\tilde{\Phi}_{1-\overline{\rho}}(x)
\right]\,.
\end{split}
\end{equation}
\end{itemize}
In equations (a) and (b), the integrals as well as the sums (or series)
with respect to the Blaschke zeroes all converge in the absolute sense.}

{\bf Proof:} We define the functions 
$f_{11}, f_{12}, f_{21}, f_{22} : (1,\infty) \to \R$ by
\begin{equation}\label{fmn}
\begin{split}
f_{11}(x) & := \frac{2 \sqrt{x}}{\pi}\int \limits_{0}^{\infty} 
\frac{ \pi N(u) - \vartheta(u)-2 \arctan(2u)}{u^2+\frac14} 
\cdot \sin(u \log x) \, du\,, \\
f_{12}(x) & := -\sum \limits_{\substack{\rho \,:\, \zeta(\rho)=0\\ \mbox{Re}(\rho) > \frac12}}
\left[
\tilde{\Phi}_{\rho}(x)+\tilde{\Phi}_{1-\overline{\rho}}(x)
-2\tilde{\Phi}_{\frac12(1+\rho-\overline{\rho}\,)}(x)
\right]\,, \\
f_{21}(x) & := -\frac{2 \sqrt{x}}{\pi}\int \limits_{0}^{\infty} 
\frac{ \log |\zeta(\frac12+iu)|}{u^2+\frac14} 
\cdot \cos(u \log x) \, du\,, \\
f_{22}(x) & := -\sum \limits_{\substack{\rho \,:\, \zeta(\rho)=0\\ \mbox{Re}(\rho) > \frac12}}
\left[
\tilde{\Phi}_{\rho}(x)-\tilde{\Phi}_{1-\overline{\rho}}(x)
\right]\,. \\
\end{split}
\end{equation}
The absolute convergence of the $f_{21}$-integral was already mentioned
for the formulation of the Balazard-Saias-Yor Theorem, 
and the absolute convergence of the $f_{11}$-integral
results from the first equation in \eqref{nt1}.
For the absolute convergence of the sums $f_{12}, f_{22}$
we notice the third equation in \eqref{nt1}
and for $x>1$, $\mbox{Im}(\alpha) \neq 0$ the asymptotic law
\eqref{phiasym1}, which implies that
\begin{equation}\label{phiasym2}
\begin{split}
 \tilde{\Phi}_{\alpha}(x) = 
\frac{1}{\alpha (\alpha -1)}\,\frac{x^{\alpha}}{\log x} 
+ R(\alpha,x)\,,\quad
 |R(\alpha,x)| \leq  \, 
\frac{2\,x^{\mbox{Re}(\alpha)}}{|\mbox{Im}(\alpha)|^2 \, {\log}^2 x}\,.
\end{split}
\end{equation} 
We also obtain that the convergence of the expressions
in \eqref{fmn} is uniform on each compact interval $0 < x_0 \leq x \leq x_1$,
such that the functions $f_{11}, f_{12}, f_{21}, f_{22}$ are continuous.

Using the following asymptotic behaviour for $x \to \infty$,
\begin{equation}\label{pnt}
\pi_*(x)-\mbox{Li}(x) = O(x\,e^{-c \sqrt{\log x}})\,, \quad 
c>0 \mbox{~constant}\,,
\end{equation}
we conclude from \eqref{zetap1} in the limit $s \to 1$ for real $s>1$ that
\begin{equation}\label{fsternasym}
\lim \limits_{x \to \infty} \frac{f_{*}(x)}{x}=
\int \limits_1^{\infty} \frac{\pi_*(t)-\mbox{Ei}_0(\log t)}{t^2}\,dt = 0\,.
\end{equation}
In accordance with this we also obtain from \eqref{phiasym2} that
\begin{equation}\label{fmnasymptotik}
\lim \limits_{x \to \infty} \frac{f_{mn}(x)}{x}=0\,, \quad m, n \in \{1,2\}\,.
\end{equation}
On the other hand we have from the definition of $f_*$ 
\begin{equation}\label{fsternlim1}
\lim \limits_{x \to 1} f_{*}(x)=0\,,
\end{equation}
and the symmetry of the zeta zeroes with respect
to complex conjugation and \eqref{bsyintegral},
\eqref{phialpha}, \eqref{phialphamod} imply that
\begin{equation}\label{fmnlim1}
\begin{split}
\lim \limits_{x \to 1} f_{11}(x)=0\,, \quad 
\lim \limits_{x \to 1} f_{12}(x)= 0\,,\\
\lim \limits_{x \to 1} f_{21}(x)= -2 \Omega_{\zeta}\,, \quad 
\lim \limits_{x \to 1} f_{22}(x)=  2 \Omega_{\zeta}\,.\\
\end{split}
\end{equation}
From Theorem (2.2), \eqref{zetap2} and \eqref{Phimellin},
\eqref{phialphamod} we obtain for $\mbox{Re}(s)>1$ 
with a constant $K \in \C$ that
\begin{equation}\label{hilf1}
1=\exp \left[\, s(s-1)\int \limits_{1}^{\infty}
\frac{f_*(x)-\left(f_{11}(x)+f_{12}(x) + K \cdot x  \right)}{x^{s+1}}\,dx \,
\right]\,.
\end{equation}
We conclude for $x>1$ from \eqref{fsternasym} to \eqref{fmnlim1}
with Mellin's inversion formula
\begin{equation}\label{hilf2}
f_*(x)=f_{11}(x)+f_{12}(x)\,.
\end{equation}
Thus we have shown the first part of Theorem (3.3).\\

In order to prove the second part, we first recall 
that the function $\zeta_{B}$
in \eqref{burnfact0} and \eqref{zetabpoisson}
is analytic for $\mbox{Re}(s)>\frac{1}{2}$ with
\begin{equation}\label{expzb1}
\begin{split}
B(1) \zeta_{B}(s) & = \frac{s-1}{s}\,\zeta(s) \, \frac{B(1)}{B(s)}\\
 & = \exp \left[ -\Omega_{\zeta} +
\frac{2}{\pi} \left(s-\frac12\right)
\int \limits_{0}^{\infty}
\frac{\log \left| \zeta(\frac12+iu) \right|} 
{u^2+(s-\frac12)^2} \, du\,
\right] \,.\\
\end{split}
\end{equation}
Since we have the limit
\begin{equation}\label{expzb1}
B(1) \lim \limits_{s \to 1}
\zeta_{B}(s) = \frac{B(1)}{B(1)}=1\,,
\end{equation}
the real valued function $Q : (1/2, \infty) \to \R$ with
\begin{equation}\label{expzb2}
Q(s)  := \frac{\log \left( B(1)\zeta_{B}(s) \right)}{s(s-1)}\\
\end{equation}
has an analytical continuation $Q_*$ in $\mbox{Re}(s)> \frac12$ with
\begin{equation}\label{expzb3}
\begin{split}
Q_*(s) & = \frac{2\Omega_{\zeta}}{s} +
\frac{2}{\pi}\frac{s-\frac12}{s(s-1)}
\int \limits_0 ^{\infty} \left[ 
\frac{\log \left| \zeta(\frac12+iu) \right|}{(u^2+(s-\frac12)^2)} -
\frac{\log \left| \zeta(\frac12+iu) \right|}{(u^2+\frac14)} 
\right] \, du\\
& = \frac{2\Omega_{\zeta}}{s}-\frac{2}{\pi}(s-\frac12)
\int \limits_0 ^{\infty} \frac{\log \left| \zeta(\frac12+iu) \right|} 
{(u^2+\frac14) \, (u^2+(s-\frac12)^2)} \, du\\
& =  \int \limits_1^{\infty} 
\frac{2 \Omega_{\zeta} + f_{21}(x)}{x^{s+1}}\,dx\,.\\
\end{split}
\end{equation}
For the last equation in \eqref{expzb3} 
we have used Fubini's theorem with
\begin{equation}\label{expzb4}
\int \limits_1^{\infty} 
\frac{\sqrt{x} \cos(u \log x)}{x^{s+1}}\,dx = 
\frac{s-\frac12}{u^2+(s-\frac12)^2}\,, \quad \mbox{Re}(s)>\frac12\,,
\end{equation}
because all integrals converge in the absolute sense.

Next we obtain from \eqref{Phimellin}
and \eqref{burnfact2} for $\mbox{Re}(s)>1$
that
\begin{equation}\label{Bsdurchb1}
\frac{B(s)}{B(1)}=\exp
\left[
s(s-1) \int \limits_1 ^{\infty}
\{
2 \Omega_{\zeta}\cdot (x-1) -
\sum \limits_{\substack{\rho \,:\, \zeta(\rho)=0\\ \mbox{Re}(\rho) > \frac12}}
\left(
\Phi_{\rho}(x)-\Phi_{1-\overline{\rho}}(x)
\right)\,\}
\frac{dx}{x^{s+1}}
 \right]\,,
\end{equation}
and with \eqref{zetap2}, \eqref{burnfact0}, \eqref{expzb2}, \eqref{expzb3}
for $\mbox{Re}(s)>1$
\begin{equation}\label{finale1}
\begin{split}
\frac{s-1}{s}\,\zeta(s) &=
\exp
\left[
s(s-1) \int \limits_1 ^{\infty}
\frac{f_*(x)}{x^{s+1}}\,dx
 \right]\\
&= \left( \zeta_B(s) B(1) \right) \cdot \frac{B(s)}{B(1)}\\
&=\exp
\left[
s(s-1) \int \limits_1 ^{\infty}
\frac{2 \Omega_{\zeta}+f_{21}(x)}{x^{s+1}}\,dx
 \right] \cdot \frac{B(s)}{B(1)}\,.
\end{split}
\end{equation}
Note the relation \eqref{phialphamod} between 
$\Phi_{\alpha}(x)$ and 
$\tilde{\Phi}_{\alpha}=x \cdot \varphi_{\alpha -1}(x)-\varphi_{\alpha}(x)$. 
Then we combine \eqref{Bsdurchb1} and \eqref{finale1}
and use \eqref{fsternasym}-\eqref{fmnlim1} to conclude
the second part of the theorem from Mellin's inversion formula.
\dokend\\

Finally we mention without going into details
that the mathematical technique developed
here can be used as well in order to derive an integral form for
other explicite formulas related to the distribution of prime
numbers, for example for the following result\\

{\bf Theorem (3.4)}~{\em For $x>1$, $\alpha,r \in \C$ and $u \in \R$ we define} 

\begin{equation}\label{pir}
\begin{split}
{\Theta}(x,{\alpha}) &:= 
\left\{\begin{displaystyle}
\begin{array}{cl}  
- \D \frac{x^{\alpha}-1}{\alpha}+\mbox{Ei\,}_0(\alpha \log x) \log x
\,, & \alpha \neq 0\,, \\
-\log x\,, & \alpha = 0\,,\\
\end{array}\end{displaystyle}
\right.\\
K(x,r,u) &:= \left \{
\begin{array}{cl}  
\frac{1}{ \pi}\,\left[
\frac{x^{1/2+iu-r}-1}{(1/2+iu-r)^2}-\frac{\log x}{1/2+iu-r}
\right]\,, & r \neq 1/2+iu\,,\\
\frac{\log ^2 x}{2 \pi}\,, & r = 1/2+iu\,,
\end{array}\right.\\
\pi_{*,r}(x) &:= \sum \limits_{n \leq x} \frac{\Lambda(n)}{n^r \log n}\,,
\quad
\psi_{r}(x) := \sum \limits_{n \leq x} \frac{\Lambda(n)}{n^r}\,.\\
\end{split}
\end{equation}
{\em Then we have}
\begin{equation}\label{pir}
\begin{split}
\int \limits_{1}^{x} \frac{\pi_{*,r}(y)}{y}\,dy
&=\pi_{*,r}(x) \log x - \psi_r(x)\\
&=\Theta(x,1-r)-\Theta(x,-r)\\
&+ \int \limits_{0}^{\infty} 
\left( K(x,r,u) +  K(x,r,-u) \right)
\log |\zeta(\frac12+iu)|\,du\\
&-\sum \limits_{\substack{\rho \,:\, \zeta(\rho)=0\\ \mbox{Re}(\rho) > \frac12}}
\left\{
\Theta(x,\rho-r)-\Theta(x,1-\overline{\rho}-r)
\right\}\,.
\end{split}
\end{equation} 
{\em On the right hand side in \eqref{pir} the integral as well as the sum 
(or series) with respect to the Blaschke zeroes 
converge in the absolute sense.} \dokend


\end{document}